\documentclass[12pt]{article}

\input{amssym.def}
\input{amssym}
\usepackage{eufrak}

\newtheorem{theorem}{Theorem}
\newtheorem{corollary}[theorem]{Corollary}
\newtheorem{definition}[theorem]{Definition}
\newtheorem{example}[theorem]{Example}
\newtheorem{proposition}[theorem]{Proposition}
\newtheorem{remark}[theorem]{Remark}

\def\be{\begin{equation}}
\def\ee{\end{equation}}

\def\ot{\otimes}
\def\vv{V^{\ot 2}}
\def\C{{\Bbb C}}
\def\R{{\Bbb R}}
\def\K{{\Bbb K}}
\def\Sym{{\rm Sym\, }}
\def\qq{q^{-1}}
\def\End{{\rm End\, }}
\def\De{\Delta}
\def\Tr{{\rm Tr}}
\def\Ren{R_{\End(V)}}
\def\h{{\hbar}}
\def\lrqh{{\cal L}(R_q,\h)}
\def\lrq{{\cal L}(R_q)}
\def\lrqo{{\cal L}(R_q,1)}
\def\gq{\gggg_{q}}
\def\gggg{\frak g}
\def\span{{\rm span\, }}

\begin{document}

\title{Quantum Lie algebras via modified Reflection Equation
Algebra}
\author{
\rule{0pt}{7mm} Dimitri
Gurevich\thanks{gurevich@univ-valenciennes.fr}\\
{\small\it USTV, Universit\'e de Valenciennes,
59304 Valenciennes, France}\\
\rule{0pt}{7mm} Pavel Saponov\thanks{Pavel.Saponov@ihep.ru}\\
{\small\it Division of Theoretical Physics, IHEP, 142284
Protvino,
Russia} }

\date{}

\maketitle

\section{Introduction}
\label{sec:int}

A Lie super-algebra was historically the first generalization of
the notion of a Lie algebra. Lie super-algebras were introduced
by
physicists in studying dynamical models with fermions. In
contrast
with the usual Lie algebras defined via the classical flip  $P$
interchanging any two elements $P(X\ot Y)=Y\ot X$, the definition
of a Lie super-algebra is essentially based on a super-analog of
the permutation $P$. This super-analog is defined on a
$\rm{Z}_2$-graded vector space $V=V_{\overline 0}\oplus
V_{\overline 1}$ where ${\overline 0},{\overline 1}\in \rm{Z}_2$
is a "parity". On homogeneous elements (i.e. those belonging to
either $V_{\overline 0}$ or $V_{\overline 1}$) its action is
$P(X\ot Y)=(-1)^{{\overline X}{\overline Y}} Y\ot X$, where
${\overline X}$ stands for the parity of a homogeneous element
$X\in V$.

Then a Lie super-algebra is the following data
$$
(\gggg=\gggg_{\overline 0}\oplus \gggg_{\overline 1},
\,P:\gggg\ot\gggg\to
\gggg\ot\gggg,\, [\,\,,\,]:\gggg\ot \gggg \to \gggg)\,,
$$
where $\gggg$ is a super-space, $P$ is a super-flip, and
$[\,\,,\,]$
is a Lie super-bracket, i.e. a linear operator which is subject
to
three axioms:
\begin{enumerate}
\item $[X,Y]=-(-1)^{{\overline X}{\overline Y}}[Y,X]$; \item
$[X,[Y,Z]]+(-1)^{{\overline X}({\overline Y}+ {\overline
Z})}[Y,[Z,X]] +(-1)^{{\overline Z} ({\overline X}+{\overline
Y})}[Z,[X,Y]]=0$; \item ${\overline {[X,Y]}}={\overline
X}+{\overline Y}$.
\end{enumerate}
Here $X,Y,Z$ are assumed to be arbitrary  homogenous elements of
$\gggg$. Note that all axioms can be rewritten via the
corresponding
super-flip. For instance the axiom 3 takes the form
$$P(X\ot [Y,Z])= [\,\,,\,]_{12}  P_{23}P_{12}(X\ot Y\ot Z).$$
(As usual, the indices indicate the space(s) where a given
operator is applied.)

In this paper we discuss the problem what is a possible
generalization of the notion of a Lie super-algebra related to
"flips" of more general type.

The first generalization of the notion of a Lie super-algebra was
related to gradings different from $\rm{Z}_2$. The corresponding
Lie type algebras were called $\Gamma$-graded  ones (cf.
\cite{Sh}).

The next step was done in \cite{G1} where there was introduced a
new generalization of the Lie algebra notion related to an
involutive symmetry defined as follows. Let $V$ be a vector space
over a ground field $\K$ (usually $\C$ or $\R$) and $R:\vv\to\vv$
be a linear operator. It is called a {\em braiding} if it
satisfies the  quantum Yang-Baxter equation
$$
R_{12}R_{23}R_{12}=R_{23}R_{12}R_{23}\,,
$$
where $R_{12}=R\ot I,\,\,R_{23}=I\ot R$ are operators in the
space
$V^{\ot 3}$. If such a braiding satisfies the condition $R^2=I$
(resp.,  $(R-q\, I)(R+q^{-1}\,I)=0,\,q\in\K)$  we call it an {\em
involutive symmetry} (resp., a {\em Hecke symmetry}). In the
latter case $q$ is assumed to be generic\footnote{Note that there
exists a big family of Hecke and involutive symmetries which are
not deformations of the usual flip (cf. \cite{G2}). Even the
Poincar\'e-Hilbert (PH) series corresponding to the "symmetric"
$\Sym(V)=T(V)/\langle \Im(q I-R)\rangle$ and "skew-symmetric"
$\bigwedge(V)=T(V)/\langle \Im(\qq I+R)\rangle$ algebras  can
drastically  differ from the classical ones, whereas the PH
series
are stable under a deformation.}.

Two basic examples of  {\em generalized} Lie algebras  are
analogs
of the Lie algebras $gl(n)$ and $sl(n)$ (or of their
super-analogs
$gl(m|n)$ and $sl(m|n)$). They can be associated to any
"skew-invertible" (see Section 2) involutive symmetry
$R:\vv\to\vv$. We denote them $gl(V_R)$ and $sl(V_R)$
respectively. The generalized Lie algebras $gl(V_R)$ and
$sl(V_R)$
are defined in the space $\End(V)$ of endomorphisms of the space
$V$. Their enveloping algebras $U(gl(V_R))$ and $U(sl(V_R))$
(which can be defined in a natural way) are equipped with a
braided Hopf structure such that the coproduct coming in its
definition acts on the generators $X\in gl(V_R)$ or $sl(V_R)$ in
the classical manner: $\De:X\to X\ot 1+1\ot X$.

Moreover, if an involutive symmetry $R$ is a deformation of the
usual flip (or super-flip) the enveloping algebras  $U(gl(V_R))$
and $U(sl(V_R))$ are deformations of their classical (or super-)
counterparts.

There are known numerous attempts to define a quantum (braided)
Lie algebra similar to generalized ones but without assuming $R$
to be involutive. Let us mention some of them: \cite{W},
\cite{LS}, \cite{DGG}, \cite{GM}. In this note we compare the
objects defined there  with $gl$ type Lie algebras-like objects
introduced recently in \cite{GPS}. Note that the latter objects
can be associated with any  skew-invertible Hecke symmetry, in
particular, that related to Quantum Groups (QG) of $A_n$ series.
Their enveloping algebras are treated in terms of the modified
reflection equation algebra (mREA) defined bellow. These
enveloping algebras have good deformation properties and the
categories of their finite dimensional {\em equivariant}
representations look like those of the Lie algebras $gl(m|n)$.
Moreover, these  algebras  can be equipped with a structure of
braided bi-algebras. Though the corresponding coproduct acts on
the generators of the algebras in a non-classical way it is in a
sense intrinsic (it has nothing in common with the coproduct in
the QGs). Moreover, it allows to define braided analogs of
(co)adjoint vectors fields.

We think that apart from generalized Lie algebras related to
involutive symmetries (described in Section 2) there is no
general
definition of a quantum (braided) Lie algebra. Moreover,
reasonable  quantum Lie algebras exist only for the $A_n$ series
(or more generally, for any skew-invertible Hecke symmetry). As
for the quantum Lie algebras of the $B_n,\, C_n,\, D_n$ series
introduced in \cite{DGG}, their enveloping algebras are not
deformations of their classical counterparts and for this reason
they are somewhat pointless objects.
\vskip 3mm

{\noindent\bf Acknowledgement} The work of P.Saponov was
partially supported by RFBR grant no. 05-01-01086.

\section{Generalized Lie algebras}
\label{sec:2}

Let $R:\vv\to\vv$ be an involutive symmetry. Then the data
$$
(V, R, [\,\,,\,]:\vv \to V)
$$
is called a {\em generalized Lie algebra} if the following holds
\begin{enumerate}
\item $[\,\,,\,]\,R(X\ot Y)=-[X,Y]$; \item
$[\,\,,\,]\,[\,\,,\,]_{12}(I+R_{12}R_{23}+R_{23}R_{12}) (X\ot
Y\ot
Z)=0$; \item $R [\,\,,\,]_{12} (X\ot Y \ot Z)= [\,\,,\,]_{12}
R_{23}R_{12}(X\ot Y\ot Z)$.
\end{enumerate}
Such a generalized Lie algebra is denoted $\gggg$.

Note that the generalized Jacobi identity (the axiom 2) can be
rewritten in one of the following equivalent forms
\begin{itemize}
\item[]
$[\,\,,\,]\,[\,\,,\,]_{23}(I+R_{12}R_{23}+R_{23}R_{12})(X\ot Y\ot
Z)=0$; \item[] $[\,\,,\,]\,[\,\,,\,]_{12}(X(Y\ot Z-R(Y\ot
Z)))=[X,[Y,Z]]$;
\item[] $[\,\,,\,]\,[\,\,,\,]_{23}((X\ot Y-R(X\ot
Y))Z)=[[X,Y],Z]$.
\end{itemize}

\begin{example}\rm If $R$ is the the usual flip
then the third axiom is fulfilled automatically and we get a
usual
Lie algebra. If $R$ is a super-flip then we get a Lie
super-algebra. In the both cases $R$ is involutive.
\end{example}

The enveloping algebras of the generalized Lie algebra $\gggg$
can
be defined in a natural way:
$$
U(\gggg)=T(V)/\langle X\ot Y-R(X\ot Y)-[X,Y]\rangle\,.
$$
(Hereafter $\langle I\rangle$ stands for the ideal generated by a
set $I$.) Let us introduce the symmetric algebra $\Sym(\gggg)$ of
the generalized Lie algebra $\gggg$ by the same formula but with
0
instead of the bracket in the denominator of the above formula.

For this algebra there exists a version of the
Poincar\'e-Birhoff-Witt theorem.

\begin{theorem}
 The algebra $U(\gggg)$ is canonically isomorphic to
$Sym(\gggg)$.
\end{theorem}

A proof can be obtained via the Koszul property established in
\cite{G2} and the results of \cite{PP}. Also, note that similarly
to the classical case  this isomorphism can be realized via a
symmetric (w.r.t. the symmetry $R$) basis.

\begin{definition}
\label{def:3} We say  that a given braiding $R:\vv\to\vv$ is
skew-invertible if there exists a morphism $\Psi:\vv\to\vv$ such
that
$$
\Tr_2\Psi_{12} R_{23}=P_{13}=\Tr_2\Psi_{23} R_{12}
$$
where $P$ is the usual flip.
\end{definition}

If $R$ is a skew-invertible braiding, a "categorical
significance"
can be given to the dual space of $V$. Let $V^*$ be the vector
space dual to $V$. This means that there exist a non-degenerated
pairing $\langle\,\,,\,\,\rangle : V^*\ot V\to \K$ and an
extension of the symmetry $R$ to the space $(V^*\oplus V)^{\ot
2}\to (V^*\oplus V)^{\ot 2}$ (we keep the same notation for the
extended braiding) such that the above pairing is $R$-invariant.
This means that on the space $V^*\ot V\ot W$ (resp., $W\ot V^*\ot
V$)  where either $W=V$ or $W=V^*$ the following relations hold
$$
R\,\langle\,\,,\,\,\rangle_{12}=\langle\,\,,\,\,\rangle_{23}\,
R_{12}\,R_{23}\qquad({\rm resp.,}\quad
R\,\langle\,\,,\,\,\rangle_{23}=\langle\,\,,\,\,\rangle_{12}\,
R_{23}\,R_{12})\,.
$$
(Here as usual, we identify $X\in W$ with $X\ot 1$ and $1\ot X$.)

Note that if such an extension exists it is unique. By fixing
bases $x_i\in V$ and $x_i\ot x_j\in \vv$ we can identify the
operators $R$ and $\Psi$ with matrices $\|R_{ij}^{kl}\|$ and
$\|\Psi_{ij}^{kl}\|$ respectively. For example,
$$R(x_i\ot x_j)=R_{ij}^{kl}\,x_k\ot x_l$$
(from now on we assume the summation over the repeated indices).

Then the above definition can be presented in the following
matrix
form
$$R_{ij}^{kl}\, \Psi_{lm}^{\,jn}=\delta_m^k\delta_i^n.$$
If $ ^ix$ is the {\em left dual} basis of the space $V^*$, i.e.
such that  $\langle\, ^jx,x_i\rangle=\delta_i^j$ then we put
$$\langle x_i,\,
^jx\rangle=\langle\,\,,\,\,\rangle\,\Psi_{ik}^{jl}\,\, { ^kx}\ot
x_l= C_i^j,\quad {\rm where} \quad C_i^j=\Psi_{ik}^{jk}.$$ (Note
that the operator $\Psi$ is a part of  the braiding $R$ extended
to the space $(V^*\oplus V)^{\ot 2}$.) By doing so, we ensure
$R$-invariance of the pairing $V\ot V^*\to \K$.

As shown in \cite{GPS} for any skew-invertible Hecke symmetry $R$
the following holds
$$ C_i^j\, B_j^k=q^{-2a}\delta_i^k,\quad {\rm where} \quad
B_i^j=\Psi_{ki}^{kj}$$ with an integer $a$ depending on the the
HP
series of the algebra $\Sym(V)$ (see footnote 1). So, if
$q\not=0$
the operators $C$ and $B$ (represented by the matrices $
\|C_i^j\|$ and  $\|B_j^k\|$ respectively) are invertible.
Therefore, we get a non-trivial pairing
$$\langle\,\,,\,\,\rangle: (V\oplus V^*)^{\ot 2}\to\K$$
which is $R$-invariant.

Note that these operators $B$ and $C$ can be introduced without
fixing any basis in the space $V$ as follows
 \be B_2 = Tr_{(1)}(\Psi_{12}),\quad C_1 = Tr_{(2)}(\Psi_{12}).
 \label{B-mat} \ee
Let us exhibit  an evident but very important property of these
operators
 \be Tr_{(1)}(B_1\,R_{12}) = I,\quad Tr_{(2)}(C_2\, R_{12}) = I.
 \label{tr-BR} \ee

By fixing the basis $h_i^j=x_i\ot { ^jx}$ in the space
$\End(V)\cong V\otimes V^*$  equipped with the usual product
$$
\circ:\End(V)^{\ot 2}\to\End(V)
$$
we get the following multiplication table $h_i^j\circ
h_k^l=\delta_k^j\, h_i^l$.

Below we use another basis in this algebra, namely that
$l_i^j=x_i\ot x^j$ where $x^j$ is the {\em right dual} basis in
the space $V^*$, i.e. such that
 $\langle x_i, x^j \rangle=\delta_i^j$.
Note that  the multiplication table for the the product  $\circ$
in this basis is $l_i^j\circ l_k^m=B_k^j\, l_i^m$ (also see
formula (\ref{l-action})).

Let $R$ be the above extension of a skew-invertible braiding to
the space $(V^*\oplus V)^{\ot 2}$. Then a braiding $\Ren :
\End(V)^{\ot 2}\to\End(V)^{\ot 2}$ can be defined in a natural
way:
$$
\Ren=R_{23}\,R_{34}R_{12}\,R_{23}\,,
$$
where we used the isomorphism $\End(V)\cong V\ot V^*$.

Observe that the  product $\circ$ in the space $\End(V)$ is
$R$-invariant and therefore $\Ren$-invariant. Namely, we have
$$
\Ren(X\circ Y, Z)=\circ_{23} (\Ren)_{12} (\Ren)_{23} (X\ot Y  \ot
Z)\,,
$$
$$
\Ren(X,Y \circ  Z)=\circ_{12} (\Ren)_{23} (\Ren)_{12} (X\ot Y \ot
Z)\,.
$$

\begin{example}
\label{ex:4} \rm Let $R:\vv\to \vv$ be a skew-invertible
involutive symmetry. Define a generalized Lie bracket by the rule
$$
[X,Y]=X\circ Y-\circ \Ren (X\ot Y)\,.
$$
Then  the data $(\End(V), \Ren, [\,\,,\,])$ is a generalized Lie
algebra (denoted  $gl(V_R)$).

Besides, define the $R$-trace $\Tr_R: \End(V)\to \K$ as follows
$$
\Tr_R(h_i^j)=B_j^i\,h_i^j ,\qquad X\in \End(V)\,.
$$

The $R$-trace possesses the following properties :

\begin{itemize}
\item The pairing
$$\End(V)\ot \End(V)\to\K:\,\,X\otimes Y\mapsto
\langle X,Y\rangle=Tr_R (X\circ Y)$$
 is non-degenerated;

\item It is $\Ren$-invariant in the following sense
$$
\Ren((\Tr_R X)\ot Y)=(I\ot \Tr_R)\Ren(X\ot Y)\,,
$$
$$
\Ren (X\ot (\Tr_R Y))=(\Tr_R\ot I)\Ren (X\ot Y)\,;
$$

\item $\Tr_R \,[\,\,,\,]=0$.
\end{itemize}

Therefore the set $\{X\in gl(V_R)| Tr_R \, X=0\}$ is closed
w.r.t.
the above bracket. Moreover,  this subspace squared is invariant
w.r.t the symmetry $\Ren$. Therefore this subspace (denoted
$sl(V_R)$) is a generalized Lie subalgebra.
\end{example}

Observe that the enveloping algebra of any generalized Lie
algebra
possesses a braided Hopf algebra structure such that the
coproduct
$\Delta$ and antipode $S$ are defined on the generators in the
classical way
$$
\De(X)=X\ot 1+1\ot X,\qquad S(X)=-X\,.
$$
For details the reader is referred to \cite{G2}.

Also, observe that while $R$ is  a super-flip the generalized Lie
algebra $gl(V_R)$ (resp., $sl(V_R)$) is nothing but the Lie
super-algebras $gl(m|n)$ (resp., $sl(m|n)$).

\section{Quantum Lie algebras for $B_n,\, C_n, \, D_n$ series}
\label{sec:3}

In this Section we restrict ourselves to the braidings coming
from
the QG $U_q(\gggg)$ where $\gggg$ is a Lie algebra of one of the
series $B_n,\, C_n, \, D_n$. By  the Jacobi identity, the usual
Lie bracket
$$
[\,\,,\,]: \gggg\ot \gggg\to \gggg
$$
is a $\gggg$-morphism.

Let us equip the space $\gggg$  with a $U_q(\gggg)$ action which
is a
deformation of the usual adjoint one. The space $\gggg$ equipped
with such an action is denoted $\gq$. Our immediate goal is to
define an operator
$$
[\,\,,\,]_q:\gq\ot \gq\to \gq
$$
which would be a $U_q(\gggg)$-covariant deformation of the
initial
Lie bracket. This means that the $q$-bracket satisfies the
relation
$$
[\,\,,\,]_q (a_1(X)\ot a_2(Y))=a([X\ot Y]_q)\,,
$$
where $a$ is an arbitrary element of the QG $U_q(\gggg)$, $a_1\ot
a_2=\De(a)$ is the Sweedler notation for the QG coproduct $\De$,
and $a(X)$ stands for the result of applying the element $a\in
U_q(\gggg)$ to an element $X\in \gq$.

Let us show that $U_q(\gggg)$-covariance of the bracket entails
its
$R$-invariance where $R=P\,\pi_{\gggg\ot\gggg}({\cal R})$ is the
image
of the universal quantum $R$-matrix ${\cal R}$ composed with the
flip $P$. Indeed, due to the relation
$$
\De_{12}({\cal R})={\cal R}_{13}{\cal R}_{23}\,,
$$
we have (by omitting $\pi_{\gggg\ot\gggg}$)
\begin{eqnarray*}
&& R[\,\,,\,]_1(X\ot Y\ot Z)=P{\cal R} ([X,Y]\ot Z) =
P[\,\,,\,]_{12}\De_{12}{\cal R}(X\ot Y\ot Z)=\\
&& P[\,\,,\,]_{12} {\cal R}_{13}{\cal R}_{23}(X\ot Y\ot Z)
= P[\,\,,\,]_{12}P_{13}R_{13}P_{23}R_{23} (X\ot Y\ot Z)=\\
&& P[\,\,,\,]_{12}P_{13}P_{23}R_{12}R_{23}(X\ot Y\ot
Z)=[\,\,,\,]_{23}R_{12}R_{23}(X\ot Y\ot Z)\,.
\end{eqnarray*}

Finally, we have
$$R[\,\,,\,]_{12}=[\,\,,\,]_{23}R_{12}R_{23},\,\,\,R[\,\,,\,]_{23
}=[\,\,,\,]_{12}R_{23}R_{12}$$ (the second relation can be
obtained in a similar way).

Thus, the $U_q(\gggg)$-covariance  of the bracket $[\,\,,\,]_q$
can
be considered as an analog of the axiom 3 from the above list. In
fact, if $\gggg$ belongs to one of the series $B_n,\, C_n$ or
$D_n$,
this property suffices for unique (up to a factor) definition of
the bracket $[\,\,,\,]_q$. Indeed, in this case it is known  that
if one extends the adjoint action of $\gggg$ to the space
$\gggg\otimes \gggg$ (via the coproduct in the enveloping
algebra),
then the latter space is multiplicity free with respect to this
action. This means that there is no isomorphic   irreducible
$\gggg$-modules in the space $\gggg\ot \gggg$. In particular, the
component isomorphic to $\gggg$ itself appears only in the
skew-symmetric subspace of $\gggg\ot \gggg$. A similar property
is
valid for decomposition of the space $\gq\ot \gq$ into a direct
sum of irreducible $U_q(\gggg)$-modules (recall that $q$ is
assumed
to be generic).

Thus, the map $[\,\,,\,]_q$, being a $U_q(\gggg)$-morphism, must
kill all components in the decomposition of $\gggg_q\ot \gggg_q$
into
a direct sum of irreducible $\gggg_q$-submodules except for the
component isomorphic to $\gq$. Being restricted to this
component,
the map $[\,\,,\,]_q$ is an isomorphism. This property uniquely
defines the map $[\,\,,\,]_q$ (up to a non-zero factor). For an
explicit computation of the structure constants of the
$q$-bracket
$[\,\,,\,]_q$ the reader is referred to the paper \cite{DGG}.
Note
that the authors of that paper embedded the space $\gq$ in the QG
$U_q(\gggg)$. Nevertheless, it is possible to do all the
calculations without such an embedding but using the QG just as a
substitute of the corresponding symmetry group.

Now, we want to define  the enveloping algebra of a quantum Lie
algebra $\gq$. Since the space $\gq\ot \gq$ is multiplicity free,
we conclude that there exists a unique $U_q(\gggg)$-morphism
$P_q:\gq\ot \gq\to \gq\ot \gq$ which is a deformation of the
usual
flip and such that $P_q^2=I$. Indeed, in order to introduce such
an operator it suffices to define q-analogs of symmetric and
skew-symmetric components in $\gq\ot \gq$. Each of them can be
defined as a direct sum of irreducible $U_q(\gggg)$-submodules of
$\gq\ot \gq$ which are $q$-counterparts of the
$U_q(\gggg)$-modules
entering  the usual symmetric and skew-symmetric subspaces
respectively.

Now, the enveloping algebra can be defined as a quotient
$$
U(\gq)=T(\gq)/\langle X\ot Y-P_q(X\ot Y)-[\,\,,\,]_q\rangle\,.
$$

Thus, we have defined the quantum Lie algebra $\gq$ and its
enveloping algebra related to the QG of $B_n,\, C_n, \, D_n$
series. However, the question what properties of these quantum
Lie
algebras are similar to those of generalized Lie algebras is
somewhat pointless since the algebra $U(\gq)$ is not a
deformation
of its classical counterpart. Moreover, its "$q$-commutative"
analog (which is defined similarly to the above quotient but
without the $q$-bracket $[\,\,,\,]_q$ in the denominator) is not
a
deformation of the algebra $\Sym(\gggg)$. For the proof, it
suffices
to verify that the corresponding semiclassical term is not a
Poisson bracket. (However, it becomes Poisson bracket upon
restriction to the corresponding algebraic group.)

\begin{remark} \rm A similar construction of a quantum Lie
algebra is valid for any skew-inver\-tible braiding of the
Birman-Murakami-Wenzl type. But for the same reason it is out of
our interest.
\end{remark}

Also, note that the Lie algebra $sl(2)$ possesses a property
similar to that above: the space $sl(2)\ot sl(2)$ being equipped
with the extended adjoint action is a multiplicity free
$sl(2)$-module. So, the corresponding quantum Lie algebra and its
enveloping algebra can be constructed via the same scheme.
However, the latter algebra is a deformation of its classical
counterpart. This case is consider in  the next Sections as a
part
of our general construction related to Hecke symmetries.

\section{Modified Reflection Equation Algebra and its
representation theory}

In this section we shortly describe the modified reflection
equation algebra (mREA) and the quasitensor Schur-Weyl category
of
its finite dimensional {\em equivariant} representations. Our
presentation is based on the work \cite{GPS}, where these objects
were considered in full detail.

The starting point of all constructions is a Hecke symmetry $R$.
As was mentioned in Introduction, the Hecke symmetry is a linear
operator $R:\vv\to\vv$, satisfying the quantum Yang-Baxter
equation and the additional Hecke condition
$$
(R-q\,I)(R+q^{-1}I) = 0\,,
$$
where a nonzero $q\in {\Bbb K}$ is generic, in particular, is not
a primitive root of unity. Besides, we assume $R$ to be
skew-invertible (see Definition \ref{def:3}).

Fixing bases $x_i\in V$ and $x_i\ot x_j\in \vv$, $1\le i,j\le
N=\dim V$, we identify $R$ with a $N^2\times N^2$ matrix
$\|R_{ij}^{kl}\|$. Namely, we have
 \be R  (x_{i}\ot x_{j})= R_{ij}^{kl}\,
x_{k}\ot x_{l} \,, \label{flip} \ee where the lower indices label
the rows of the matrix, the upper ones --- the columns.

As is known, the Hecke symmetry $R$ allows to define a
representations $\rho_R$ of the $A_{k-1}$ series Hecke algebras
$H_k(q)$, $k\ge 2$, in tensor powers $V^{\otimes k}$:
$$\rho_R: H_k(q)\rightarrow {\rm End}(V^{\otimes k})\quad
\rho_R(\sigma_i) = R_i:=I^{\otimes (i-1)}\otimes R\otimes
I^{\otimes(k-i-1)}\,,
$$
where elements $\sigma_i$, $1\le i\le k-1$ form the set of the
standard generators of $H_k(q)$.

The Hecke algebra $H_k(q)$ possesses the primitive idempotents
$e^\lambda_a\in H_k(q)$, which are in one-to-one correspondence
with the set of all standard Young tableaux $(\lambda,a)$,
corresponding to all possible partitions $\lambda\vdash k$. The
index $a$ labels the tableaux of a given partition $\lambda$ in
accordance with some ordering.

Under the representation $\rho_R$, the primitive idempotents
$e_a^\lambda$ are mapped into the projection operators \be
E^{\lambda}_a(R)=\rho_R(e^\lambda_a)\in \End(V^{\otimes k})\,,
\label{Y-proj} \ee these projectors being some polynomials in
$R_i$, $1\le i\le k-1$.

Under the action of these projectors the spaces $V^{\otimes k}$,
$k\ge 2$, are expanded into the direct sum \be V^{\otimes k} =
\bigoplus_{\lambda\vdash k} \bigoplus_ {a=1}^{d_\lambda}
V_{(\lambda,a)}, \qquad V_{(\lambda,a)} = {\rm
Im}(E_a^\lambda)\,,
\label{V-decom} \ee where the number $d_\lambda$ stands for the
total number of the standard Young tableaux, which can be
constructed for a given partition $\lambda$.

Since the projectors $E^\lambda_a$ with different $a$ are
connected by invertible transformations, all spaces
$V_{(\lambda,a)}$ with fixed $\lambda$ and different $a$ are
isomorphic. Note, that the isomorphic spaces $V_{(\lambda,a)}$
(at
a fixed $\lambda$) in decomposition (\ref{V-decom}) are treated
as
particular embeddings of the space $V_\lambda$ into the tensor
product $V^{\otimes k}$. Hereafter we use the notation
$V_\lambda$
for the class of the spaces $V_{(\lambda,a)}$ equipped with one
or
another embedding in $V^{\otimes k}$.

In a similar way we define  classes $V^*_\mu$. First, note that
the Hecke symmetry being extended to the space $(V^*)^{\ot 2}$ is
given in the basis $x^{i}\ot x^j$ as follows
$$R(x^{i}\ot x^j)=R^{ji}_{lk}\,\,x^{k}\ot x^l$$
(and similarly in the basis $ ^ix \ot \, ^jx$). It is not
difficult to see that the operator $R$ so defined in the space
$(V^*)^{\ot 2}$ is a Hecke symmetry. Thus, by using the above
method we can introduce  spaces $V_{(\mu,a)}^*$ looking like
those
from (\ref{V-decom}) and define the classes $V^*_\mu$.

Now, let us define a rigid quasitensor Schur-Weyl category ${\rm
SW}(V)$ whose objects are spaces $V_\lambda$ and $V^*_\mu$
labelled by partitions of  nonnegative integers, as well as their
tensor products $V_\lambda\ot V^*_\mu$ and all finite  sums of
these spaces.

Among the morphisms of the category ${\rm SW}(V)$ are the above
left and right pairings and the set of braidings
$R_{U,W}:U\otimes
W\to W\otimes U$ for any pair of objects $U$ and $W$. These
braidings can be defined in a natural way. In order to define
them
on a couple of objects of  the form $V_\lambda\ot V^*_\mu$
 we embed them into appropriate products
$V^{\ot k}\ot (V^*)^{\ot l}$ and define the braiding $R_{U,W}$ as
an appropriate restriction. Note, that all these braidings are
$R$-invariant maps (cf. \cite{GPS} for detail). Note that the
category ${\rm SW}(V)$ is monoidal quasitensor rigid according to
the standard terminology (cf. \cite{CP}).

Now we are aiming  at  introducing modified reflection equation
algebra and  equipping objects of the category ${\rm SW}(V)$ with
a structure of its modules.

Again, consider the space $\End(V)$ equipped with the basis
$l_i^{\,j}$ (see Section 2). Note that  the element $l_i^{\,j}$
acts on the elements of the space $V$ as follows
 \be
l_i^{\,j}(x_k):=x_i\,\langle x^j,x_k\rangle = x_iB_k^{\,j}\,.
\label{l-action} \ee Introduce the $N\times N$ matrix $L =
\|l_i^{\,j}\|$. Also, define its ''copies'' by the iterative rule
\be L_{\overline 1}:= L_1:=L\otimes I,\qquad L_{\overline{k+1}}:=
R_kL_{\overline k}R_k^{-1}\,. \label{iterat} \ee Observe that the
isolated spaces $L_{\overline{k}}$  have no meaning (except for
that $L_{\overline{1}}$). They can be only correctly understood
in
the products $L_{\overline{1}}L_{\overline{2}}$,
$L_{\overline{1}}L_{\overline{2}}L_{\overline{3}}$ and so on, but
this notation is useful in what follows.

\begin{definition}
\label{def:10} \rm The associative algebra generated by the unit
element $e_{\cal L}$ and the indeterminates $l_i^{j}$ $1\le
i,j\le
N$ subject to the following matrix relation
 \be R_{12} L_{1}R_{12} L_{1} - L_{1}R_{12}L_{1}R_{12} -
\hbar\,(R_{12}\,L_{1} - L_{1}\,R_{12}) = 0\,, \label{mREA} \ee is
called the modified reflection equation algebra (mREA) and
denoted
$\lrqh$.
\end{definition}
Note, that at $\hbar = 0$ the above algebra is known as the
reflection equation algebra $\lrq$. Actually, at $q\not = \pm 1$
one has $\lrqh\cong \lrq$. Since at $\hbar\not=0$ it is always
possible to renormalize generators $L\mapsto \hbar\,L$. So, below
we
consider  the case $\hbar =1$.

Thus, the mREA is the quotient algebra of the free tensor algebra
$T(\End(V))$ over the two-sided ideal, generated by the matrix
elements of the left hand side of (\ref{mREA}). It can be shown,
that the relations (\ref{mREA}) are $R$-invariant, that is the
above two-sided ideal is invariant when commuting with any object
$U$ under the action of the braidings $R_{U,\End(V)}$ or
$R_{\End(V),U}$ of the category ${\rm SW(V)}$.

Taking into account (\ref{tr-BR}) one can easily prove, that the
action (\ref{l-action}) gives a basic (vector) representation of
the mREA $\lrqo$ in the space $V$ \be
\rho_1(l_i^{\,j})\triangleright x_k = x_iB_k^{\,j}\,,
\label{l-rep} \ee where the symbol $\triangleright$ stands for
the
(left) action of a linear operator onto an element. Since $B$ is
non-degenerated, the representation is irreducible.

Another basic (covector) representation $\rho_1^*:{\cal L}(R_q,1)
\rightarrow \End(V^*)$ is given by \be
\rho_1^*(l_i^j)\triangleright x^k = -x^r\, R_{ri}^{kj}\,.
\label{duals} \ee one can prove, that the maps $\End(V)\to
\End(V)$ and $\End(V)\to \End(V^*)$ generated by
$$
l_i^{\, j}\mapsto\rho_1(l_i^{\, j}) \quad{\rm and}\quad l_i^{\,
j}\mapsto\rho^*_1(l_i^{\, j})
$$
are the morphisms of the category ${\rm SW}(V)$.

\begin{definition} A representation $\rho:{\cal L}(R_q,1)
\rightarrow \End(U)$ where $U$ is an object of the category ${\rm
SW}(V)$ is called equivariant if its restriction to $\End(V)$ is
a
categorical morphism.
\end{definition}

Thus, the above representations $\rho_1$ and $\rho_1^*$ are
equivariant.

Note that there are known representations of the mREA which are
nor equivariant. However, the class of equivariant
representations
of the mREA is very important. In particular, because the tensor
product of two equivariant $\lrqo$-modules  can be also  equipped
with a structure of an equivariant $\lrqo$-module via a "braided
bialgebra structure" of the mREA.

Let us briefly describe this structure. It consists of two maps:
the braided coproduct $\Delta$ and counit $\varepsilon$.

The coproduct $\Delta$ is an algebra homomorphism of $\lrqo$ into
the associative algebra $\mbox{\bf L}(R_q)$ which is defined as
follows.

\begin{itemize}
\item As a vector space over the field $\Bbb K$ the algebra
$\mbox{\bf L}(R_q)$ is isomorphic to the tensor product of two
copies of mREA
$$
\mbox{\bf L}(R_q) = {\cal L}(R_q,1) \otimes {\cal L}(R_q,1)\,.
$$
\item The product $\star : (\mbox{\bf L}(R_q))^{\otimes
2}\rightarrow \mbox{\bf L}(R_q)$ is defined by the rule \be
(a_1\otimes b_1)\star (a_2\otimes b_2):=a_1 a^\prime_2 \otimes
b^\prime_1 b_2\,,\qquad a_i\otimes b_i \in \mbox{\bf L}(R_q)\,,
\label{br-pr} \ee where $a_1a^\prime_2$ and $b_1b^\prime_2$ are
the usual product of mREA elements, while $a^\prime_1$ and
$b^\prime_1$ result from the action of the braiding $R_{{\rm
End}(V)}$ (see Section \ref{sec:2}) on the tensor product
$b_1\otimes a_2$ \be a^\prime_2\otimes b^\prime_1:= R_{{\rm
End}(V)} (b_1\otimes a_2)\,. \label{ti-el} \ee
\end{itemize}

The braided coproduct $\Delta$ is now defined as a linear map
$\Delta: {\cal L}(R_q,1)\rightarrow \mbox{\bf L}(R_q)$ with the
following properties: \be
\begin{array}{l}
\Delta(e_{\cal L}):= e_{\cal L}\otimes e_{\cal L}\\
\rule{0pt}{6mm} \Delta(l_i^j) := l_i^j\otimes e_{\cal L}+ e_{\cal
L}\otimes l_i^j -
(q-q^{-1})\sum_k l_i^k\otimes l_k^j\\
\rule{0pt}{6mm}
\Delta(ab):=\Delta(a)\star\Delta(b)\qquad\forall\,a,b\in {\cal
L}(R_q,1)\,.
\end{array}
\label{copr} \ee

In addition to (\ref{copr}), we introduce a linear map
$\varepsilon :{\cal L}(R_q,1) \rightarrow {\Bbb K}$ \be
\begin{array}{l}
\varepsilon(e_{\cal L}):= 1\\
\rule{0pt}{6mm}
\varepsilon(l_i^j):= 0\\
\rule{0pt}{6mm} \varepsilon(ab):= \varepsilon(a)\varepsilon(b)
\qquad \forall\,a,b\in {\cal L}(R_q,1)\,.
\end{array}
\label{coed} \ee

One can show (cf. \cite{GPS}) that the maps $\Delta$ and
$\varepsilon$ are indeed algebra homomorphisms and that they
satisfy the relation
$$
({\rm id}\otimes \varepsilon)\,\Delta = {\rm id} =
(\varepsilon\otimes {\rm id})\, \Delta\,.
$$

Let now $U$ and $W$ be two equivariant mREA-modules with
representations $\rho_U:{\cal L}(R_q,1)\rightarrow \End(U)$ and
$\rho_W:{\cal L}(R_q,1)\rightarrow \End(W)$ respectively.
Consider
the map $\rho_{U\otimes W}:\mbox{\bf L$(R_q)$} \rightarrow
\End(U\otimes W)$ defined as follows \be \rho_{U\otimes
W}(a\otimes b)\triangleright (u\otimes w) =
(\rho_U(a)\triangleright u^\prime)\otimes (\rho_W(b^\prime)
\triangleright w)\,, \qquad a\otimes b\in \mbox{\bf L$(R_q)$}\,,
\label{rep-ll} \ee where
$$
u^\prime \otimes b^\prime:=R_{\End(V),U}(b\otimes u)\,.
$$
Definition (\ref{rep-ll}) is self-consistent since the map
$b\mapsto \rho_W(b^\prime)$ is also a representation of the mREA
${\cal L}(R_q,1)$.

The following proposition holds true.

\begin{proposition} {(\rm \cite{GPS})}
\label{pro:29} The action {\rm (\ref{rep-ll})} defines a
representation of the algebra $\mbox{\bf L}(R_q)$.
\end{proposition}

Note again, that  the equivariance of the representations in
question plays a decisive role in the proof of the above
proposition.

As an immediate corollary of the proposition \ref{pro:29} we get
the rule of tensor multiplication of equivariant $\lrqo$-modules.

\begin{corollary}
Let $U$ and $W$ be two ${\cal L}(R_q,1)$-modules with equivariant
representations $\rho_U$ and $\rho_W$. Then the map ${\cal
L}(R_q,1)\rightarrow \End(U\otimes W)$ given by the rule  \be
{a\mapsto \rho_{U\otimes W}(\Delta(a))\,, \qquad \forall\,a\in
{\cal L}(R_q,1)}\, {\rm \label{rep-prod}} \ee is an equivariant
representation. Here the coproduct $\Delta$ and the map
$\rho_{U\otimes W}$ are given respectively by formulae {\rm
(\ref{copr})} and {\rm (\ref{rep-ll})}.
\end{corollary}

Thus, by using (\ref{rep-prod}) we can extend the basic
representations $\rho_1$ and $\rho_1^*$ to the representations
$\rho_k$ and $\rho_l^*$ in tensor products $V^{\otimes k}$ and
$(V^{*\otimes l})$ respectively. These representations are
reducible, and their restrictions on the representations
$\rho_{\lambda,a}$ in the invariant subspaces $V_{(\lambda,a)}$
(see (\ref{V-decom})) are given by the projections \be
\rho_{\lambda,a} = E^\lambda_a\circ\rho_k\, \label{restr} \ee and
similarly for the subspaces $V^*_{(\mu,a)}$. By using
(\ref{rep-prod}) once more we can equip each object of the
category  ${\rm SW}(V)$ with the structure of an equivariant
${\cal L}(R_q,1)$-module.

\section{Quantum Lie algebras related to Hecke symmetries}

In this section we consider the question to which extent one can
use the scheme of section \ref{sec:2} in the case of
non-involutive Hecke symmetry $R$ for definition of the
corresponding  Lie algebra-like object. For such an object
related
to a Hecke symmetry $R$ we use the term {\it quantum} or {\it
braided} Lie algebra. Besides, we require the mREA, connected
with
the same symmetry $R$, to be an analog of the enveloping algebra
of the quantum Lie algebra. Finally, we compare the properties of
the above generalized Lie algebras and quantum ones.

Let us recall the interrelation of a usual  Lie algebra $\frak g$
and its universal enveloping algebra $U(\frak g)$. As is known,
the universal enveloping algebra for a Lie algebra $\frak g$ is a
unital associative algebra $U(\frak g)$ possessing the following
properties:
\begin{itemize}
\item There exists a linear map $\tau:\frak g\to U(\frak g)$ such
that 1 and ${\rm Im}\,\tau$ generate the whole $U(\frak g)$.

\item The Lie bracket $[x,y]$ of any two elements of $\frak g$
has
the image
$$
\tau([x,y]) = \tau(x)\tau(y) - \tau(y)\tau(x).
$$
\end{itemize}

Let us rewrite these formulae in an equivalent form. Note that
the
tensor square $\frak g\otimes \frak g$ splits into the direct sum
of symmetric and skew symmetric components
$$
\frak g\otimes \frak g = \frak g_s\oplus \frak g _a, \qquad \frak
g_s = {\rm Im}\,{\cal S},\quad \frak g_a  = {\rm Im}\, {\cal
A}\,,
$$
where ${\cal S}$ and ${\cal A}$ are the standard
(skew)symmetrizing operators
$$
{\cal S}(x\otimes y) = x\otimes y + y\otimes x, \quad {\cal
A}(x\otimes y) = x\otimes y - y\otimes x\,,
$$
where we neglect the usual normalizing factor $1/2$. Then the
skew-symmetry property of the classical Lie bracket is equivalent
to the requirement \be [\,\,,\,]\,{\cal S}(x\otimes y) = 0\,.
\label{sk-sy} \ee

The image of the bracket in $U(\frak g)$ is presented as follows
\be \tau([x,y]) = \circ {\cal A}\,(\tau(x)\otimes \tau(y))\,,
\label{im-br} \ee where $\circ$ stands for the product in the
associative algebra $U(\frak g)$.
\begin{itemize}
\item The Jacobi identity for the Lie bracket $[\,\,,\,]$
translates into the requirement that the correspondence $x\mapsto
[x,\,]$ generate the (adjoint) representation of $U(\frak g)$ in
the linear subspace $\tau(\frak g)\subset U(\frak g)$.
\end{itemize}

So, we define a braided Lie algebra as a linear subspace ${\cal
L}_1=\End(V)$ of the mREA $\lrqo$, which generates the whole
algebra
and is equipped with the quantum Lie bracket. We want the bracket
to
satisfy some skew-symmetry condition, generalizing (\ref{sk-sy}),
and  define a representation of the mREA in the same linear
subspace ${\cal L}_1$ via an analog of the Jacobi identity.

As ${\cal L}_1$, let us take the linear span of mREA generators
$$
{\cal L}_1 = \End(V)\cong V\otimes V^*\,.
$$
Together with the unit element this subspace generate the whole
$\lrqo$ by definition.

In order to find the quantum Lie bracket, consider a particular
representation of $\lrqo$ in the space $\End(V)$. In this case
the
general formula (\ref{rep-prod}) reads
$$
l_i^{\,j}\mapsto \rho_{V\otimes V^*}(\Delta(l_i^{\,j}))\,,
$$
where  we should take the basic representations (\ref{l-rep}) and
(\ref{duals}) as $\rho_{V}(l_i^j)$ and $\rho_{V^*}(l_i^j)$
respectively. Omitting straightforward calculations, we write the
final result in the compact matrix form \be \rho_{V\otimes
V^*}(L_{\overline 1})\triangleright L_{\overline 2} = L_1R_{12} -
R_{12}L_1\,, \label{ad} \ee where the matrix $L_{\overline k}$ is
defined in (\ref{iterat}).

Let us define \be [L_{\overline 1},L_{\overline 2}] = L_1R_{12} -
R_{12}L_1\,. \label{q-lie-br} \ee The generalized skew-symmetry
(the axiom 1 from Section \ref{sec:2}) of this bracket is now
modified as follows. In the space ${\cal L}_1\otimes{\cal L}_1$
one can construct two projection operators ${\cal S}_q$ and
${\cal
A}_q$ which are interpreted as $q$-symmetrizer and
$q$-skew-symmetrizer respectively (cf. \cite{GPS}). Then
straightforward calculations show that the above bracket
satisfies
the relation \be [\,\,,\,]\,{\cal S}_q(L_{\overline 1}\otimes
L_{\overline 2}) = 0\,, \label{q-skew} \ee which is the
generalized skew-symmetry condition, analogous to (\ref{sk-sy}).

Moreover, if we rewrite the defining commutation relations of the
mREA
(\ref{mREA}) in the equivalent form \be L_{\overline
1}L_{\overline 2} - R_{12}^{-1} L_{\overline 1}L_{\overline 2}
R_{12} = L_1R_{12} - R_{12}L_1\,, \label{eq-form} \ee we come to
a
generalization of the formula (\ref{im-br}).

By introducing an operator
$$Q : {\cal L}_1^{\ot 2}\to {\cal L}_1^{\ot 2},\quad
Q(L_{\overline 1}L_{\overline 2})= R^{-1}_{12}L_{\overline
1}L_{\overline 2}R_{12}
$$
we can present the relation (\ref{eq-form}) as follows \be
L_{\overline 1}L_{\overline 2}-Q (L_{\overline 1}L_{\overline
2})=[L_{\overline 1},L_{\overline 2}]. \label{Q-form} \ee It
looks
like the defining relation of the enveloping algebra of a
generalized Lie algebra. (Though we prefer to use the notations
$L_{\overline k}$ it is possible to exhibit the maps $Q$ and
$[\,\,,\,\,]$ in the basis $l_i^j\ot l_k^m$.) Observe that the
map
$Q$ is a braiding. Also, note that the operators ${\cal S}_q$ and
${\cal A}_q$ can be expressed in terms of $Q$ and its inverse
(cf.
\cite{GPS}).

We call the data $(\gggg={\cal L}_1, Q, [\,\,,\,\,])$ the $gl$
type
quantum (braided) Lie algebra. Note that if $q=1$ (i.e. the
symmetry $R$ is involutive) then $Q=\Ren$ and this quantum Lie
algebra is nothing but the generalized Lie algebra $gl(V_R)$ and
the corresponding mREA becomes isomorphic to its enveloping
algebra.

Let us list the properties of the the quantum Lie algebra in
question.

\begin{itemize}
\item The bracket $[\,\,,\,\,]$ is skew-symmetric in the sense of
(\ref{q-skew}). \item The q-Jacobi identity is valid in the
following form \be [\,\,,\,][\,\,,\,]_{12} =
[\,\,,\,][\,\,,\,]_{23}(I-Q_{12})\,. \label{q-Jac} \ee \item The
bracket $[\,\,,\,\,]$ is $R$-invariant. Essentially, this means
that the following relations hold
$$
R_{\End(V)}[\,\,,\,]_{23} =
[\,\,,\,]_{12}(R_{\End(V)})_{23}(R_{\End(V)})_{12},$$ \be
R_{\End(V)}[\,\,,\,]_{12} =
[\,\,,\,]_{23}(R_{\End(V)})_{12}(R_{\End(V)})_{23}\,.
 \label{last}
\ee
\end{itemize}

So, the adjoint action
$$
L_{\overline 1}\triangleright L_{\overline 2}=[L_{\overline 1},
L_{\overline 2}]
$$
is indeed a representation. By chance (!) the representation
$\rho_{V\ot V^*}$ coincides with this adjoint action.

Turn now to the question of the "$sl$-reduction", that is, the
passing from the mREA ${\cal L}(R_q,1)$ to the quotient algebra
\be
{\cal SL}(R_q):= {\cal L}(R_q,1)/\langle \Tr_RL\rangle\,, \qquad
\Tr_RL:= \Tr(C L)\,, \label{sl-quo} \ee (see Section 2 for the
operator $C$). The element $\ell:=\Tr_RL$ is central in the mREA,
which can be easily proved by calculating the $R$-trace in the
second space of the matrix relation (\ref{mREA}).

To describe the quotient algebra ${\cal SL}(R_q)$ explicitly, we
pass to a new set of  generators $\{f_i^j,\ell\}$, connected with
the initial one by a linear transformation: \be l_i^j = f_i^j +
(\Tr(C) )^{-1}\delta_i^j\,\ell\quad {\rm or}\quad L=F+(\Tr(C)
)^{-1}I\,\ell\,, \label{shift} \ee where $F=\|f_i^j\|$. Hereafter
we assume that $\Tr\, C=\ell_i^i\not=0$. (So, the Lie
super-algebras
$gl(m|m)$ and their q-deformations are forbidden.) Obviously,
$\Tr_RF = 0$, i.e. the generators $f_i^j$ are dependent.

In terms of the new generators, the commutation relations of the
mREA
read
$$
\left\{
\begin{array}{l}
\displaystyle R_{12}F_1 R_{12}F_1 -  F_1 R_{12}F_1 R_{12} =
(e_{\cal L} - \frac{\omega}{Tr(C)}\,\ell)(R_{12}F_1
- F_1 R_{12})\\
\rule{0pt}{6mm} \ell \,F = F\,\ell\,,
\end{array}
\right.
$$
where $\omega = q-q^{-1}$. Now, it is easy to describe the
quotient (\ref{sl-quo}) . The matrix $F=\|f_i^j\|$ of ${\cal
SL}(R_q)$ generators satisfy the same commutation relations
(\ref{mREA}) as the matrix $L$ \be R_{12}F_1 R_{12}F_1 -  F_1
R_{12}F_1 R_{12}= R_{12}F_1 - F_1 R_{12}\,, \label{sl-rea} \ee
but
the generators $f_i^j$ are linearly dependent.

Rewriting this relation in the form similar to (\ref{Q-form}) we
can introduce an $sl$-type bracket. However for such a bracket
the
$q$-Jacobi identity fails. This is due to the fact the element
$\ell$ comes in the relations for $f_i^j$ (at $q=1$ this effect
disappears ). Nevertheless, we can construct a representation
$$
\rho_{V\ot V^*}:{\cal SL}(R_q)\to\End(V\ot V^*)
$$
which is an analog of the adjoint representation. In order to do
so, we rewrite the representation (\ref{ad}) in terms of the
generators $f_i^j$ and $\ell$. Taking  relation (\ref{shift})
into
account, we find, after a short calculation
\begin{eqnarray}
&&\rho_{V\otimes V^*}(\ell)\triangleright \ell = 0,\qquad
\rho_{V\otimes V^*}(F_1)\triangleright \ell = 0\,,\nonumber\\
&&\rho_{V\otimes V^*}(\ell)\triangleright F_1 =
-\omega\,\Tr(C)\,F_1\nonumber\\
&&\rho_{V\otimes V^*}(F_{\overline 1})\triangleright F_{\overline
2} = F_1 R_{12} -  R_{12}F_1 + \omega  R_{12}F_1 R_{12}^{-1}\,.
\label{sl-ad}
\end{eqnarray}

Namely, the last formula from this list defines the
representation
$\rho_{V\ot V^*}$. However,  in contrast with
 the mREA ${\cal L}(R_q,1)$, this map is different from that
 defined by the bracket $[\,\,,\,\,]$ reduced to the space
 $\span(f_i^j)$.
This is reason why the "$q$-adjoint" representation cannot be
presented in the form (\ref{q-Jac}). (Also, note that though
$\ell$ is central it acts in a non-trivial way on the elements
$f_i^j$.)

Moreover, any object $U$ of the category ${\rm SW(V)}$ above such
that
$$
\rho_U(\ell) = \chi \,I_U\,, \qquad \chi\in{\Bbb K}\,
$$
is a scalar operator, can be equipped with an ${\cal
SL}(R_q)$-module structure. First, let us observe that for any
representation $\rho_U:{\cal L}(R_q,1)\to\End(U)$ and for any
$z\in\K$ the map
$$
 \rho_U^{z} : {\cal L}(R_q,1) \to  \End(U),\qquad
 \rho_U^{z}(l_i^j) = z\rho_U(l_i^j) +
\delta_i^j(1-z)(q-q^{-1})^{-1}\,I_U
$$
is a representation of this algebra as well.

By using this freedom we can convert a given representation
$\rho_U:{\cal L}(R_q,1)\to\End(U)$ with the above property into
that $\rho_U^z$ such that $\rho_U^z(\ell)=0$. Thus we get a
representation of the algebra ${\cal SL}(R_q)$. Explicitly, this
representation is given by the formula \be \tilde\rho(f_i^j) =
\frac{1}{\xi}\left(\rho(l_i^j) -
(\Tr(C))^{-1}\rho(\ell)\,\delta_i^j \right)\,, \qquad \xi =
1-(q-q^{-1})(\Tr(C))^{-1}\chi\,. \label{sl-red} \ee

The data $(\span(f_i^j), Q, [\,\,,\,\,])$ where the bracket
stands
for the l.h.s. of (\ref{Q-form}) restricted to $\span(f_i^j)$ is
called the $sl$-type quantum (braided) Lie algebra.

Note that in the particular case  related to the QG $U_q(sl(n))$
this  quantum algebra can be treated in terms of \cite{LS} where
an axiomatic approach to the corresponding Lie algebra-like
object is given. However, we think that  any general axiomatic
definition of such objects is somewhat useless (unless the
corresponding symmetry is involutive). Our viewpoint is motivated
by the fact that for $B_n,\,C_n,\,D_n$ series there do not exist
"quantum Lie algebras" such that their  enveloping algebras have
good deformation properties.
 As for the $A_n$ series (or more
generally, for any skew-invertible Hecke symmetry) such objects
exist and can be explicitly  exhibited via the mREA. Their
properties differ from those listed in \cite{W, GM} in the
framework of an axiomatic approach to Lie algebra-like objects.

Completing the paper, we want to emphasize  that the above
coproduct can be
useful for definition of a "braided (co)adjoint vector field". In
the ${\cal L}(R_q,1)$ case these fields are naturally introduced
through the above adjoint action extended to the symmetric
algebra
of the space ${\cal L}_1$ by means of this coproduct. The
symmetric algebra can be defined via the above operators ${\cal
S}_q$  and ${\cal A}_q$. In the ${\cal SL}(R_q)$ case a similar
treatment is possible if $\Tr \, C\not=0$.

\end{document}